# From Operations to Elderly Care Outcomes: A Thematic Review of Industrial Engineering and Decision-Support Approaches

Shayan Farhang Pazhooh
Department of Industrial and Systems Engineering
Isfahan University of Technology
Isfahan 84156-83111, Iran
shayanfp2006mail@gmail.com | ORCID: 0009-0001-1366-7942

Fereshteh Parvaresh
Department of Industrial and Systems Engineering
Isfahan University of Technology
Isfahan 84156-83111, Iran
parvaresh@iut.ac.ir | ORCID: 0000-0001-8411-6307

***Abstract***— The rapid growth of the global aging population presents severe challenges to healthcare systems, necessitating efficient, equitable, and patient-centered care models. While Industrial Engineering and Operations Research (OR) provide robust optimization and decision-support tools to address these multidimensional complexities, current applications often remain fragmented. This paper presents a thematic review of 30 seminal studies at the intersection of OR and elderly care, categorizing the literature into home healthcare operations, polypharmacy management, and clinical chronotherapy. Our analysis highlights a significant methodological evolution from static, deterministic models toward dynamic and stochastic frameworks integrated with artificial intelligence (AI). Despite these advancements, a critical translational gap persists: the current OR literature is heavily dominated by process-level optimizations, such as staff routing, and struggles to translate these operational efficiencies into measurable clinical outcomes. Furthermore, holistic models bridging the transition between hospital and community care remain critically underexplored. To develop resilient and smart healthcare systems, this study proposes a conceptual framework that shifts the research focus from isolated operational tasks to integrated, multi-level decision-making. We emphasize the critical need for robust systems analysis, human-inclusive design, and the smartification of care through emerging digital technologies—including digital twins and large language models—to successfully bridge the gap between theoretical operational metrics and tangible patient-level health outcomes.

# I. INTRODUCTION

The global elderly population is projected to double to 2.1 billion by 2050 [1], severely straining healthcare systems. Elderly care is inherently complex due to multimorbidity, polypharmacy, and functional decline, demanding coordinated medical and social support through resilient operational frameworks. Industrial engineering offers critical tools—such as robust optimization, metaheuristics, simulation, risk analysis, and smart decision-support systems—to address these multidimensional challenges. While Industrial Engineering and Operations Research (IE/OR) applications range from home healthcare (HHC) routing and facility location to disaster evacuation and polypharmacy management, the literature remains fragmented. It often prioritizes isolated operational efficiency rather than integrated decision-making across strategic, tactical, and operational levels. This thematic review analyzes and categorizes 30 seminal articles to consolidate recent contributions in elderly care and pharmacotherapy. Specifically, we aim to (i) highlight the state of the art of IE/OR in elderly healthcare, (ii) identify methodological gaps, and (iii) propose a conceptual framework bridging operational efficiency, equity, and clinical outcomes to foster resilient, integrated care systems.

# II. METHODS

This purposive thematic review synthesizes key advancements and gaps at the intersection of industrial engineering and elderly care. Comprehensive searches across Scopus (for engineering/operations research), PubMed (for clinical literature), and Google Scholar (2005–2026) were conducted using keywords including "operations research," "elderly care," "home health care routing," "polypharmacy," and "chronotherapy." We selected 30 articles based on methodological rigor, seminal contributions, and framework relevance. This approach included general OR models from adjacent fields (e.g., chemotherapy optimization), as their core mathematical architectures are directly generalizable to geriatric care, requiring only the strategic adaptation of physiological constraints and objective functions. Articles were analyzed across five dimensions: scope (operational, tactical, strategic), methodology (e.g., modeling, reviews), data source (real-world vs. simulated), key findings, and limitations. The selected literature is structured into three thematic pillars (illustrated in Figure 1) to provide a focused overview of the field's evolution toward systems-level decision-making.

# III. FINDINGS

This section presents a detailed analysis of the reviewed literature, organized into a thematic framework. As illustrated in Figure 1, the research is divided into three main pillars: (i) home healthcare (HHC) operations, (ii) pharmacotherapy optimization, and (iii) clinical chronotherapy.

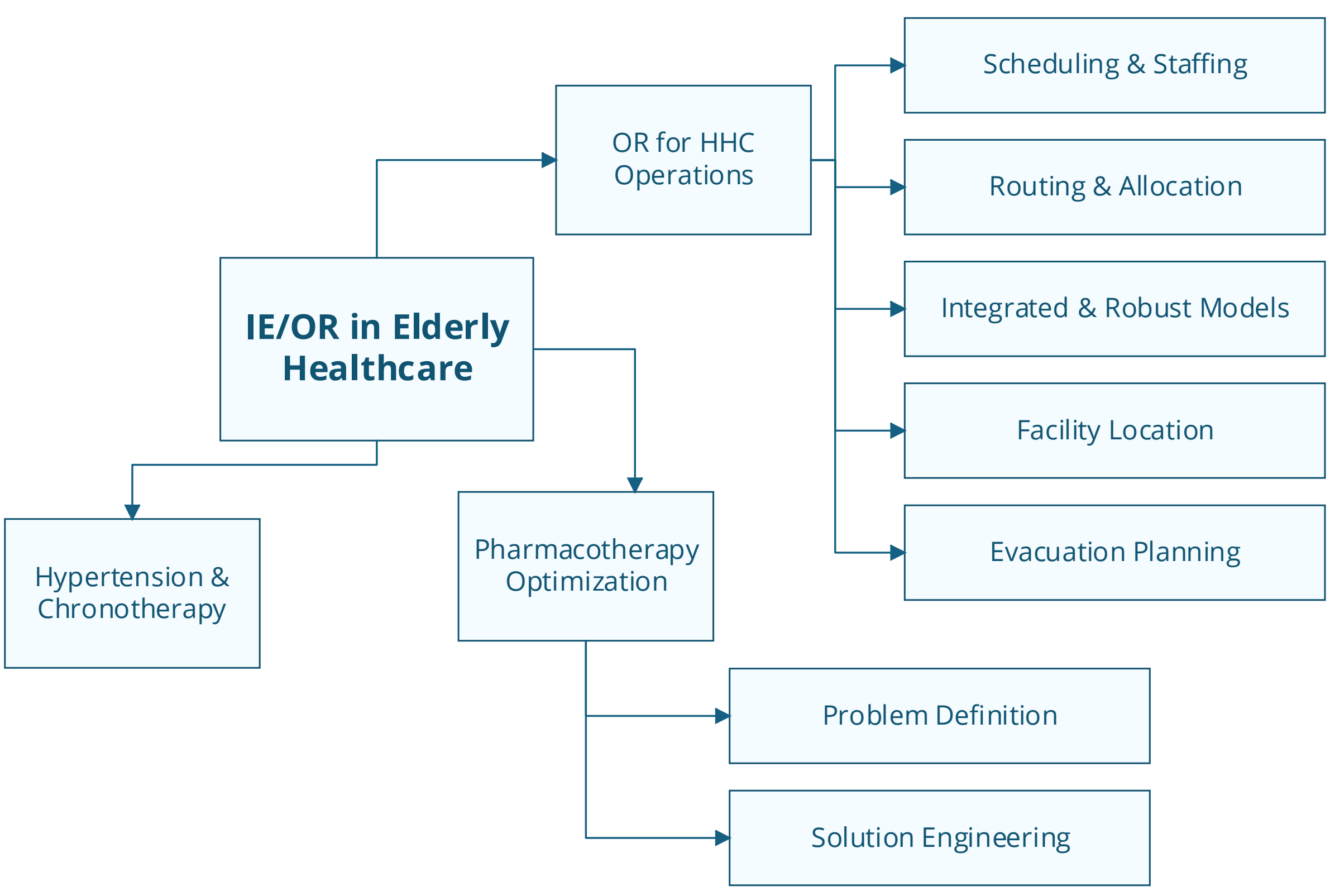

Fig. 1. Thematic Classification Framework of the Reviewed Literature

### *A. Operations Research and Optimization in Elderly Care & Home Healthcare*

As summarized in Table I, the optimization of HHC has transitioned from static, deterministic formulations toward dynamic, human-centric, and robust decision-making architectures.

The overarching trend in home healthcare operations reveals a systematic departure from deterministic, isolated tasks toward integrated, hierarchical models that actively manage uncertainty and complex human behaviors. Early scheduling literature fundamentally established the superiority of dynamic and periodic rolling-horizon policies over static assignments [4], [6], while simultaneously exposing critical trade-offs in workforce shift planning and capacity preservation [5]. As the routing domain matured, the reliance on standard metaheuristics for interdependent services [7] evolved into the deployment of sophisticated exact algorithms and hybrid architectures. For instance, managing stochastic travel and service times has been successfully addressed through two-stage integer programming and efficient L-shaped algorithms [8]. More recently, exact algorithms utilizing Benders decomposition have been introduced to seamlessly integrate staffing, routing, and scheduling under robust optimization frameworks [10].

Crucially, the optimization objectives have also shifted. Modern mathematical models no longer exclusively target cost minimization; instead, they actively incorporate complex behavioral preferences for improved social satisfaction, utilizing hybrid Markov Decision Processes (MDP) to account for patient-nurse familiarity [9]. At the strategic level, OR applications have expanded to patient-centric facility location models that prioritize minimizing travel disutility and unmet demand over mere resource utilization [11], optimize multi-period planning under tight budget constraints [12], [13], and deploy large-scale simulation-optimization frameworks for dynamic disaster evacuation planning [14]. Nevertheless, despite these significant operational and methodological strides, a structural gap persists: holistic models bridging the critical transition between hospital discharge and community-based operations remain notably underexplored [2]. This fragmentation is further corroborated by broader clustering-based systematic reviews, which expose an incoherent literature disproportionately fixated on short-term routing over long-term strategic integration [3].

TABLE I. SUMMARY OF REVIEWED LITERATURE ON OPERATIONS RESEARCH AND OPTIMIZATION

| Article (Author(s), Year) | Application Area / Problem | Methodology / Approach | Objective / Performance Measures | Key Contribution & Findings | Limitations / Future Work |
|---|---|---|---|---|---|
| Williams et al. (2021) [2] | OR in Elderly Care | Systematic Review (Scopus) | N/A | Revealing that the field is heavily fragmented and lacks holistic hospital/community integration, with an over-reliance on Markov and simulation methodologies. | 1 database (Scopus); no quality assessment. |
| Grieco et al. (2021) [3] | OR in HHC | Systematic Review (WoS), Clustering Algorithm | N/A | Revealing an incoherent literature dominated by operational routing that lacks critical strategic and tactical focus. | Keyword search limitations; subjective classification. |
| Bennett & Erera (2011) [4] | Periodic HHC Scheduling | Rolling Horizon, Heuristics | Maximize acceptance and service levels | Proved that capacity-preserving heuristics outperform immediate travel cost minimization in dynamic environments. | Single nurse; rigid appointments. |
| Håkansson (2015) [5] | Staff Scheduling | ILP, Simulation | Analyze workforce versus scheduling trade-offs | Identifying a fundamental conflict between maintaining a high full-time staff ratio and eliminating split shifts. | Assumes many available shift types. |
| Cire & Diamant (2022) [6] | Dynamic HHC Scheduling | ADP based on MDP | Maximize profit and patient acceptance | Findings prove that dynamic policies (ADP) significantly outperform static models in both metrics. | Simplified routing (TSP); lacks intraday detail. |
| Mankowska et al. (2013) [7] | Interdependent Services Routing | MILP, AVNS | Minimize travel and operational costs | Demonstrated that AVNS is highly effective for large-scale complex instances where exact solvers fail. | Deterministic times only. |
| Hashemi Doulabi et al. (2020) [8] | Synchronized & Stochastic Visits | 2-stage Stochastic Integer, L-shaped algorithm | Minimize total costs (fixed, travel, wait, penalty) | Developed the first synchronized/stochastic model, proving 2nd-stage integrality relaxation, enabling efficient L-shaped algorithm. | Implicit identical vehicle assignment. |
| Zhang et al. (2023) [9] | HHC Routing & Human Factors | Hybrid (MDP + CCP), QL-BWACO Heuristic | Maximize efficiency and social satisfaction | Integrating human factors (familiarity/preference) improved both metrics simultaneously. | Heuristic human relationships; lacks real behavioral data. |
| Naderi et al. (2023) [10] | Integrated Staffing, Routing & Scheduling | MIP, Robust Optimization, Benders decomposition (ST-LBBD) | Minimize fixed and overtime costs | Proposed a novel exact Benders decomposition algorithm, successfully quantifying both the robustness price and nurse flexibility value. | A priori generated visit patterns. |
| Johnson et al. (2005) [11] | Facility Location | Integer Programming | Minimize travel disutility | Proving that minimizing disutility/unmet demand is a | Weak demand forecasting; poor cost data. |

| | | | and unmet demand | superior objective to maximizing utilization. | |
|---|---|---|---|---|---|
| Du & Sun (2015) [12] | Multi-period HHC Location | MIP | Minimize medical cost and satisfy demand | Findings show that optimized location simultaneously increases demand satisfaction and reduces costs. | Deterministic demand; excludes staff scheduling. |
| Davari & Van Woensel (2020) [13] | Multi-period Location | MIP, ILS | Maximize demand coverage under budget constraints | Showed that ILS is highly efficient for solving large-scale multi-period location problems. | Uses generated data; lacks uncertainty. |
| Yazdani & Haghani (2023) [14] | Evacuation Planning | Flood Simulation + Multi-objective VRP | Minimize evacuation time and vehicles | Integrated sim-opt yields highly effective and dynamic large-scale evacuation plans. | Ignores shelter locations & general traffic. |

Abbreviations—ADP: Approximate Dynamic Programming, AVNS: Adaptive Variable Neighborhood Search, CCP: Chance-Constrained Programming, ILP: Integer Linear Programming, ILS: Iterated Local Search, MDP: Markov Decision Process, MILP: Mixed-Integer Linear Programming, MIP: Mixed-Integer Programming, QL-BWACO: Q-Learning based Black Widow Ant Colony Optimization, ST-LBBD: Spatio-Temporal Logic-Based Benders Decomposition, TSP: Traveling Salesperson Problem, VRP: Vehicle Routing Problem

## *B. Pharmacotherapy Optimization & Polypharmacy Management*

Table II outlines the critical domain of polypharmacy and medication management, illustrating the evolution from foundational clinical reviews toward Artificial Intelligence (AI)-augmented optimization and rigorous transitional care models.

Foundational reviews highlight a persistent gap between the well-documented harms of polypharmacy and the limited clinical efficacy of standard, single-faceted interventions [15], [16], [17], driving a necessary shift toward multifaceted "polypharmacy stewardship" [18]. To operationalize this stewardship, modern research increasingly deploys computational models—such as optimal control for multi-drug dosing [19]—and AI-driven systems (e.g., combining PM-TOM with LLMs) to actively minimize prescription risks based on STOPP/Beers criteria [20], [21]. However, translating these algorithmic optimizations into tangible endpoints remains highly challenging; safely reducing PIMs does not automatically improve primary outcomes like rehospitalization [22], and theoretically optimizing discharge prescriptions often creates a severe cognitive and adherence burden for frail patients [26]. Overcoming this translational gap demands rigorous implementation science [23] alongside systematic transitional care models, including telehealth protocols [24], the integration of general practitioners during hospitalization [25], and the strategic use of generative AI (e.g., GPT-4) to produce patient-oriented discharge summaries that significantly boost patient activation [27].

TABLE II. SUMMARY OF REVIEWED LITERATURE ON PHARMACOTHERAPY OPTIMIZATION & POLYPHARMACY MANAGEMENT

| Article (Author(s), Year) | Application Area / Problem | Methodology / Approach | Objective / Performance Measures | Key Contribution & Findings | Limitations / Future Work |
|---|---|---|---|---|---|
| Hajjar et al. (2007) [15] | Polypharmacy Clinical Review | Literature Review | N/A | Revealing a persistent disconnect between the well-documented harms of polypharmacy and the lack of proven, effective interventions to address it. | Inconsistent polypharmacy definition across literature. |
| Schlenk et al. (2008) [16] | Medication Adherence Systematic Review | Systematic Review of RCTs | N/A | Highlighting that no single intervention is superior, though tailored educational interventions combined with professional contact proved most effective. | Significant methodological flaws (selection bias, self-report). |
| Topinková et al. (2012) [17] | Pharmacotherapy Clinical Review | Narrative Review | N/A | Demonstrating mixed evidence, where pharmacist-led and behavioral interventions show promise but fail to deliver tangible clinical outcomes. | Methodological weaknesses in reviewed studies. |
| Daunt et al. (2024) [18] | Polypharmacy Conceptual Framework | Narrative Review | N/A | Concluding that single-faceted interventions have failed, thereby proposing a multifaceted, team-based polypharmacy stewardship model. | Conceptual framework; requires empirical validation. |
| Qods et al. (2025) [19] | Chemotherapy Dosing Strategy | Optimal Control, Genetic Algorithm | Minimize tumor size and delay resistance | Finding that simultaneous multi-drug administration consistently outperforms sequential strategies. | Parameters lack real individual experimental data. |
| Kulenovic et al. (2025) [20] | Polypharmacy Optimization | Case Study, AI Tool (PM-TOM + ChatGPT) | Minimize STOPP and Beers criteria risks | Demonstrated that synergistic analytical and generative AI is highly effective for optimizing complex polypharmacy. | Single case study; limited generalizability. |
| Xiao et al. (2025) [21] | Smart Medication Adherence System | System Development & Evaluation | Evaluate vs. pharmacist & general LLM | Proved that a specialized evidence-based system outperforms a generalist LLM in creating safe, personalized plans. | Unmeasured actual patient outcomes (adherence, clinical). |
| Ie et al. (2024) [22] | Hospitalized Elderly Medication Optimization | RCT | Primary: Composite (death, visits, | Successfully reduced PIMs but did not improve primary composite clinical | Single-center, open-label bias; clinician-driven. |

|  |  |  | rehospitalization). Sec: PIMs | outcomes, indicating safety without clinical efficacy. |  |
|---|---|---|---|---|---|
| Berian et al. (2024) [23] | Older Surgical Patients pCGA Implementation | Protocol (Systems Engineering & Implementation Science) | Primary: RE-AIM framework | Applies systems engineering (process mapping, user-centered design) to co-design an adaptable pCGA implementation package. | Protocol only; requires future multisite trial validation. |
| Hossain et al. (2025) [24] | Home Medication Management Televisits | Protocol (Cluster RCT) | Primary: Reduce PIMs (STOPP criteria) | Proposes a rigorous trial for a pharmacist-led televisit intervention to improve medication appropriateness. | Protocol only; veteran population limits generalizability. |
| De Guio et al. (2026) [25] | Surgical Care Transitions Polypharmacy | Prospective Multicenter Study | Maintain revised chronic treatments at 45 days | Showed that involving GPs in hospital medication reviews safely reduced inappropriate medications and maintained stability 45 days post-discharge. | Excluded highly frail patients; specific to French healthcare structure. |
| Frigaard et al. (2026) [26] | Medication Adherence Observational Cohort Study | Observational Cohort Study | Assess patient implementation of Rx-changes | Optimizing medications creates an adherence burden; notably, 61% of patients used outdated packaging with wrong instructions at home. | Small, selected group of hospitalized patients; single-center. |
| Rust et al. (2026) [27] | Patient Discharge Summaries | Single-blind RCT (using GPT-4) | Primary: Patient Activation Measure (PAM-13) | Found that LLM-generated patient-oriented discharge summaries significantly improved patient activation and were perceived as more empathetic. | Did not capture downstream clinical endpoints like readmission rates. |

Abbreviations— GP: General Practitioner, LLM: Large Language Model, PAM-13: Patient Activation Measure, pCGA: preoperative Comprehensive Geriatric Assessment, PIM: Potentially Inappropriate Medication, PM-TOM: Prescription Medication Treatment Optimization Model, RCT: Randomized Controlled Trial, RE-AIM: Reach Effectiveness Adoption Implementation Maintenance, STOPP: Screening Tool of Older Persons' Prescriptions

## *C. Hypertension and Chronotherapy Studies*

Table III summarizes the literature on hypertension chronotherapy, highlighting the tension between broad clinical guidelines and individualized therapeutic timing.

Using hypertension as a focused clinical case study, current guidelines emphasize the emergency department as a crucial juncture to safely initiate antihypertensive therapy [28]. However, the specific timing of drug administration (chronotherapy) remains contentious; while nocturnal blood pressure is a powerful risk predictor, universal chronotherapy lacks support from recent large-scale trials [29]. Nonetheless, meta-analyses of RCTs introduce critical nuance, demonstrating significant blood pressure response benefits for distinct subgroups, such as patients with high BMI [30]. To practically implement this personalized approach, modern operations research is leveraging mathematical modeling of circadian rhythms. Specifically, by utilizing ODE models fed by real-world wearable data, recent dynamic algorithms can reliably estimate a patient's circadian phase using only 5 days of smartwatch data [31], overcoming previous logistical barriers and making individualized chronotherapy highly feasible in real-world healthcare systems.

TABLE III. SUMMARY OF REVIEWED LITERATURE ON HYPERTENSION AND CHRONOTHERAPY STUDIES

| Article (Author(s), Year) | Application Area / Problem | Methodology / Approach | Objective / Performance Measures | Key Contribution & Findings | Limitations / Future Work |
|---|---|---|---|---|---|
| Roa et al. (2025) [28] | Hypertension Management Clinical Guideline | Narrative Review | N/A | Emphasizing the emergency department as a crucial juncture to safely initiate and optimize antihypertensive therapy for asymptomatic chronic hypertension. | Excludes hypertensive emergencies or severe comorbidities. |
| Parati et al. (2025) [29] | Nocturnal BP Position Paper | Position Paper / Narrative Review | N/A | Confirming nocturnal BP as a powerful risk predictor, but concluding that routine chronotherapy remains unsupported by recent large-scale trials. | Lacks conclusive data that therapeutically modifying nocturnal BP improves outcomes. |
| Kuang et al. (2025) [30] | Hypertension Chronotherapy | Meta-Analysis of RCTs | Evaluate BP response | Demonstrating that personalized chronotherapy yields significant benefits for distinct patient subgroups (e.g., high BMI, calcium channel blockers). | Limited to BP response; excludes cardiovascular event outcomes. |
| Lim et al. (2026) [31] | Chronotherapy Timing Optimization | Mathematical Modeling (ODE) + Wearable Data | Minimize data collection time for reliable circadian phase estimation | Developed a dynamic work history-based algorithm that reliably estimates circadian phases, reducing required smartwatch data from 17 to 5 days to make chronotherapy practical. | Requires continuous wearable data collection. |

Abbreviations— BMI: Body Mass Index, BP: Blood Pressure, ED: Emergency Department, ODE: Ordinary Differential Equation, RCT: Randomized Controlled Trial

# IV. DISCUSSION

Our synthesis reveals a central tension defining the field, conceptualized in Figure 2. The research landscape concentrates on two distinct dimensions: the highly developed domain of Operational Metrics & Process Efficiency (where OR models excel) and the ultimate goal of Patient-Level Clinical Outcomes. This disconnect represents a fundamental research gap. The "bridge" in our framework visualizes the key research directions—clinical validation, patient-centered design, and integrated systems—required to translate process improvements into tangible health enhancements.

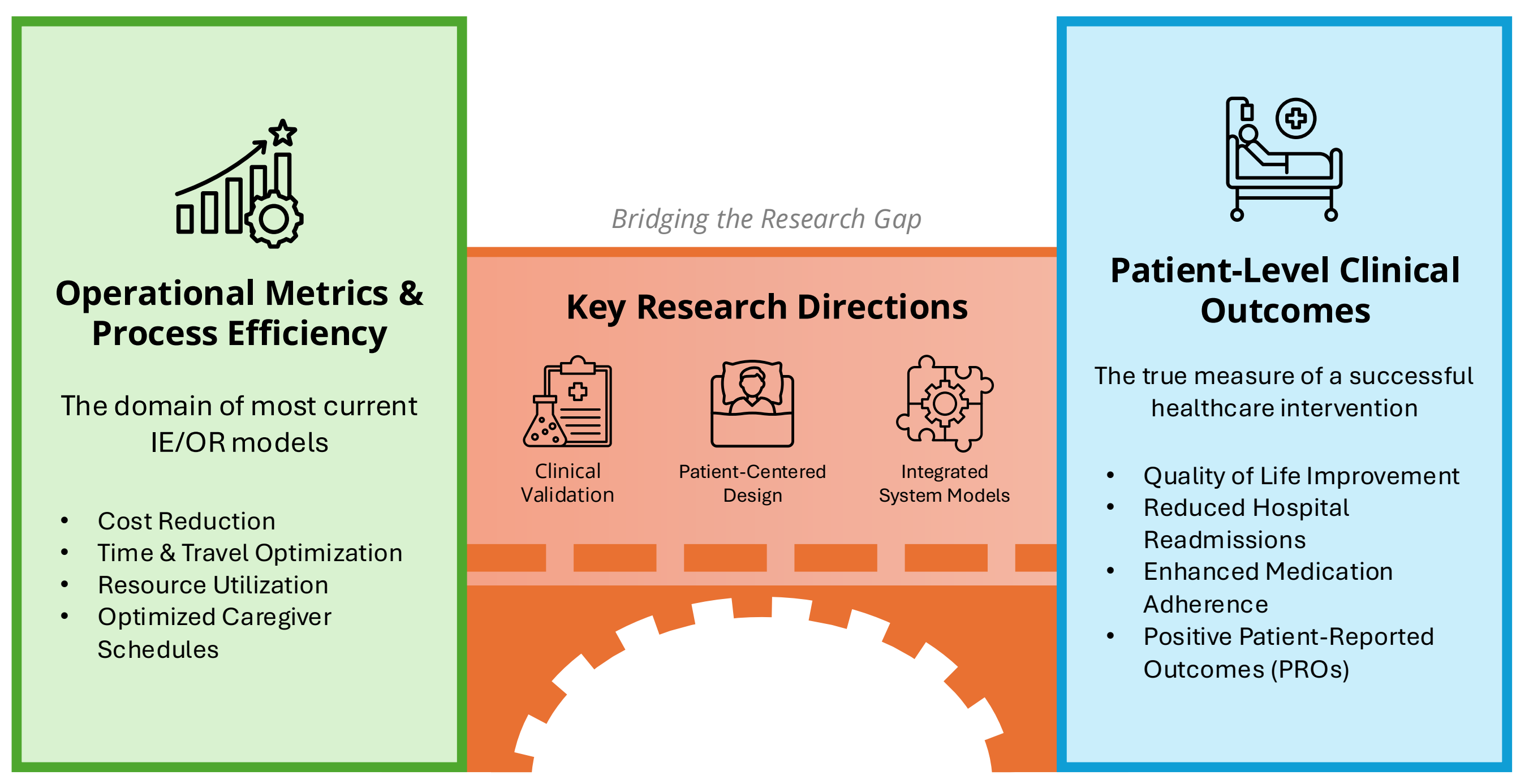

Fig. 2. Conceptual framework illustrating the research gap between operational efficiency and patient-level clinical outcomes, highlighting key research directions to bridge these two domains.

## *A. Current Trends and Persistent Limitations*

The literature demonstrates a clear shift from early descriptive risk-identification [15], [17] toward prescriptive, AI-augmented decision-support systems, integrating tools like PM-TOM with generative AI [20] and specialized LLMs [21]. Additionally, mechanistic routing models increasingly incorporate behavioral factors [9].

Despite these advances, significant limitations persist. First, research remains fragmented, disproportionately fixated on isolated operational tasks (e.g., routing, scheduling) while neglecting strategic integration across the care continuum [2], [3]. Second, models frequently rely on deterministic assumptions, simulated data [19], or single-case validations, limiting real-world scalability, though recent robust frameworks offer viable alternatives [10]. Finally, OR is fundamentally underutilized in modeling outcomes beyond pure operational efficiency. As highlighted earlier, successfully optimizing a process—such as safely reducing potentially inappropriate medications (PIMs)—does not automatically translate into improved clinical endpoints like reduced rehospitalization [22]. Unpacking this translational gap reveals three systemic bottlenecks directly relevant to OR modeling: (i) interventions often optimize "proxy metrics" by targeting low-risk medications, which mathematically improves the process but barely impacts actual patient health; (ii) mathematical models frequently overlook "competing risks" in highly frail patients, where the theoretical benefits of optimized prescriptions are overshadowed by severe underlying diseases and baseline mortality; and (iii) these optimizations are typically clinician-driven rather than patient-centered [22]. Without integrating human factors—such as ensuring patient comprehension and active behavioral modification—even the most mathematically optimized care plans will struggle to yield tangible health outcomes. This sharply underscores that healthcare delivery requires modeling beyond deterministic processes to capture the multidimensional complexity of patient realities.

### *B. Bridging the Gap: Key Research Directions*

To develop resilient healthcare systems, future research must bridge this conceptual gap (Figure 2) through three primary avenues:

- **Integrated System Models and Robust Architectures:** Overcoming the current literature fragmentation as highlighted in recent reviews [2], [3] requires hierarchical decision-support systems that seamlessly align strategic, tactical, and operational levels. To handle inherent healthcare uncertainties, frameworks must move beyond deterministic assumptions. Embedding advanced mathematical approaches—such as stochastic programming for synchronized visits [8], robust optimization for scheduling [10], and dynamic simulation-optimization [14]—is essential to ensure operational resilience.
- **Patient-Centered Design and OR-AI Hybridization:** To balance cost-efficiency with patient preferences and equity, operational models must capture complex behavioral dimensions, thereby moving beyond purely mechanistic routing [9]. Hybridizing traditional OR with AI significantly accelerates this transition. Recent literature heavily demonstrates this potential: integrating specialized LLMs and generative AI with analytical optimization tools can effectively personalize medication schedules [21] and systematically mitigate complex polypharmacy risks based on standardized clinical criteria [20]. Beyond optimizing tactical schedules and improving transitional care communication via LLM-generated, empathetic patient discharge summaries [27], the frontier of OR in elderly care lies in the deployment of Digital Twins [32]. By bridging previously siloed biomarker ecosystems—integrating continuous biosensing, genomics, and lifestyle metrics through federated learning—digital twin technologies empower predictive, proactive interventions [32]. Furthermore, integrating these AI-driven insights with OR tools allows for the translation of broad clinical guidelines into highly individualized treatment pathways. A prime example is the mathematical modeling of circadian rhythms to tailor chronotherapy for specific physiological subgroups (e.g., patients with high BMI) [30], enabling healthcare systems to definitively shift from reactive, episodic care toward adaptive, precision prevention [32].
- **Clinical Validation and Real-World Implementation:** Translating theoretical models into clinical practice demands rigorous effectiveness-implementation trials [24]. Deploying healthcare interventions is inherently a systems engineering challenge, requiring user-centered, adaptable rollouts [23]. Finally, OR models must scale to address macro-level trends. Expanding upon large-scale simulation-optimization frameworks—like those applied in disaster evacuation planning [14]—will allow policymakers to continuously evaluate and balance the economic, social, and environmental objectives of elderly care.

## V. CONCLUSION

This thematic review structured 30 key studies into a three-pillar framework, revealing a clear evolutionary trajectory in the literature: the field of elderly care operations has progressively matured from static, deterministic models designed for discrete tasks toward highly dynamic, stochastic, and AI-driven decision-support systems. However, our analysis uncovers a critical and persistent tension defining the current state of the art: a profound disconnect between process-level optimization and patient-level clinical outcomes. While OR and industrial engineering methodologies demonstrate exceptional capability in maximizing operational efficiency—such as minimizing travel costs or balancing nursing schedules— translating these theoretical gains into measurable, downstream clinical endpoints remains a formidable translational barrier.

To overcome this challenge and design sustainable, high-reliability healthcare systems, future research must shift its focus from isolated operational silos toward integrated, multi-level architectures that seamlessly connect hospital-based strategic planning with community-based operational execution. Crucially, the field demands rigorous real-world clinical validation through effectiveness-implementation trials, alongside the foundational application of systems engineering principles to ensure adaptable, user-centered rollouts. By fully embracing this patient-centered, evidence-based approach, and leveraging the hybrid power of analytical optimization and artificial intelligence, industrial engineering can fulfill its ultimate potential in driving the digital transformation of geriatric care for a rapidly aging global population.

## Declaration of Generative AI Use

We used Google Gemini 3.1 Pro for language polishing. All outputs were reviewed and verified by the authors.